\documentclass[titlepage,12pt]{article}
\usepackage{amssymb}
\usepackage{amsfonts}
\begin{document}

{\Large \bf From Incompleteness Towards \\ \\ Completeness} \\

{\bf Elem\'{e}r E Rosinger} \\
Department of Mathematics \\
and Applied Mathematics \\
University of Pretoria \\
Pretoria \\
0002 South Africa \\
eerosinger@hotmail.com \\ \\

{\bf Abstract} \\

It is argued that G\"{o}del's incompleteness theorem should be seen as self-evident, rather
than unexpected or surprising. \\ \\

{\bf A Surprising Shock and its widespread disregard} \\

Recently has been celebrated the centenary of the birth of Kurt G\"{o}del, [N]. That occasion
also recalls the 75 years since the publication of his famous paper on incompleteness. At the
time, in the 1930s, there were still strong reverberations in Mathematics following the
paradoxes in Set Theory discovered some decades earlier around the turn of the century. Three
main directions of thought had emerged in this regard. They were the so called "logicism",
following Frege, Russell and Whitehead, "formalism" supported by Hilbert, and the somewhat
less main line but rather radical "intuitionism" strongly promoted by Brouwer. \\
Logicism did not prove to be highly popular among mathematicians, since not many were ready to
accept the idea that Logic contained all of Mathematics, thus there was nothing more in
Mathematics than could be encompassed by Logic. Intuitionism proved even less popular, due to
the most severe restrictions it imposed even on Logic, let alone on Mathematics. \\

In this way, to some extent by default, Hilbert formalism came to prominence, only to be
severely questioned by G\"{o}del's 1931 incompleteness result which, in everyday terms, states
that "Every sufficiently strong axiomatic theory which is consistent must be incomplete",
[F]. \\

At the time, and to a good extent ever since, this incompleteness theorem has been considered
by so many among mathematicians who encounter it for the first time as being most surprisingly
unexpected, and in fact, rather unsettling. As it happens, however, most of the usual
mathematicians consider that, regardless of the mentioned features of the incompleteness
theorem, they can quite safely disregard it within their everyday professional pursuit. \\
Needless to say, Hilbert himself tried assiduously to diminish, if not completely dismiss, its
obvious impact on his program of formalism. \\ \\

{\bf What is self-evident should in fact not at all surprise ...} \\

And yet, what may appear to be no less surprising, this incompleteness theorem should at a
somewhat better and deeper consideration be {\it self-evident}, instead of being seen so
surprising ... \\

However, what prevents it from being self-evident is a long time and most powerfully
entrenched, even if not quite perceived "cut" in the way thinking in Mathematics, as well as
in so many other realms, is being practiced. \\
And by its very nature, Mathematics, or for that matter, Mathematical Logic, manifests this
"cut" in a most trenchant manner. \\ \\

{\bf The ever ongoing "cut" ...} \\

So then, what is this "cut" all about ? \\

For convenience, let us consider it in the case of Mathematics. Typically, and in fact, rather
without any exception, what is considered to be Mathematics, as for instance in its form
publishable in professional journals, is one thing, while mathematical thinking which leads to
such Mathematics is seen to be quite another. There may be just about one not particularly
central exception, namely, with conjectures or open problems. However, they can hardly be
published as such and all alone in professional journals. Instead, they are supposed to be
supported by a convincing argument along the lines of what is considered to be Mathematics. \\

And what is considered to be Mathematics has to be a rigorous, logically developed system
based on earlier definitions and proven results, plus possibly new definitions, and of course,
the required new results. As for axioms, they are supposed to be clear in the given context,
or occasionally, new ones may be put forward. \\
Indeed, from this point of view what is to be Mathematics is at least as clear as it may ever
happen with anything else in science. \\

But then, most obviously, hardly any mathematician does ever do Mathematics in this way. And
to make things easier to understand, let us use the analogy of what is going on even in the
very best restaurants. \\
The dishes are prepared in the kitchen, and then taken to the eating hall where the customers
sit at the neat tables, and where the dishes are served up in front of them in a spectacularly
appetizing manner. \\
Of course, the customers are not supposed to have even a look, let alone go, into the kitchen
where things are quite inevitably really messy ... \\

A similar "cut" between the "kitchen" and the "eating hall" happens in Mathematics as well,
and does so quite inevitably. \\
And needless to say, just like with restaurants, the creative mathematical activity happens in
the "kitchen", which can be the thinking of one single mathematician, or of a group who
collaborate with one another. \\ \\

{\bf The other side of the "cut" ...} \\

Strangely enough, however, what goes on in such "kitchens" is seldom brought to the attention
of the profession. And even less so of what has gone wrong in such "kitchens". Yet that
process of "cooking" new Mathematics - or in general, science - is not only essential, but
often it can also be {\it essentially different} from what is at the end seen in the "eating
halls". And in this regard we can mention several aspects, the following three of which can be
most important, no matter how much they may happen to be overlooked. \\

{\it Self-referentiality.} In present day Mathematics, or for that matter, Mathematical Logic,
so called "circular definitions" are not quite welcome, in spite of recent research in this
direction, [B-M]. However, within mathematical thinking, that is, inside the "kitchen",
nothing stops one from indulging in self-referential concepts. After all, the famous 1903
Russell paradox is but the expression of such a thinking in terms of Set Theory. \\

{\it Instantaneity.} Here, turning to Physics may better help in understanding, [R]. It is a
basic principle of both Special and General Relativity that no physical action whatever,
including information transmission can occur faster than the speed of light, [C-K]. And yet,
in the thinking of anybody, physicists or not, one can simultaneously conceive of arbitrarily
faraway objects which may possibly act upon one another. As for Quantum Mechanics, one can
conceive of two entangled quantum particles arbitrarily far from one another, yet with a
correlation of their states, which means that in one's thinking one can instantaneously
register the fact that, given the state of one of the particles, the state of the other one is
determined. \\

{\it Unattributable sources.} A good deal of original scientific ideas appear in one's mind as
if from nowhere. Furthermore, often one cannot trace precisely their origin and process of
emergence into a constituted idea even with hindsight. \\ \\

{\bf Should Incompleteness still be such a surprise ?} \\

Let us consider G\"{o}del's incompleteness theorem as relating to the Peano axioms of natural
numbers. It implies that there are statements $S$ which, together with their negation non-$S$,
do not follow as theorems in that theory. In this way, the respective theory of natural
numbers {\it branches}\, at $S$ and non-$S$, in the sense that one can add to the Peano axioms
$S$ and obtain one theory of natural numbers, or alternatively, one can add non-$S$ and,
clearly, obtain a rather different theory. \\

The decision which of such two branches to take does of course happen in the "kitchen". And
there, any of the mentioned three phenomena can manifest itself, not mention other ones not
considered above. And clearly, the third of the above phenomena - namely, the possible
presence of unattributable sources of our ideas - can be a prime source of {\it
incompleteness}. \\

Consequently, the moment we no longer consider Mathematics as being all of it contained only
on the "eating hall" side of that long ongoing and disregarded "cut", that very moment we can
only be surprised how and why we did not hit prior to G\"{o}del upon a variety of
incompleteness results. \\ \\

{\bf Not many "cuts" outside of hard science ...} \\

It may be important to enquire why such a "cut" exists, persists, and fails to be noted, let
alone, taken into account in various sciences. \\
One immediate reason can be the following. Outside of sciences - and here one means "hard
sciences" - that "cut" hardly exists. In other words, the respective "restaurants" do not have
the sharp division of "kitchen" and "eating hall". \\

Indeed, in non-intellectual realms, such as sport for instance, such a division would in fact
be highly counterproductive. \\
As for a large variety of intellectual realms, including philosophy or "soft sciences", even
if such a "cut" may be claimed to exist, that claim is questionable. And even if that "cut"
may to some extent be present, it happens - even if often regrettably so - not to be
particularly important in fact. \\

In this way, since the practice of modern "hard science" is relatively new in human history,
it may be less surprising that the mentioned "cut" has passed insufficiently noticed, let
alone considered in its implications. \\

As for the fact that Mathematics is a couple of millennia older than modern "hard science",
and nevertheless it still does not give enough attention to that "cut", one can only speculate
about the reason possibly being in the singular - and thus, insufficiently consequential -
character of Mathematics within human experience at large, as well as within the life
experience of individual mathematicians. \\ \\

\end{document}